# CONTINUOUS-TIME MEAN-VARIANCE EFFICIENCY: THE 80% RULE

By Xun Li[1] and Xun Yu Zhou[2]

*National University of Singapore and Chinese University of Hong Kong*

This paper studies a continuous-time market where an agent, having specified an investment horizon and a targeted terminal mean return, seeks to minimize the variance of the return. The optimal portfolio of such a problem is called mean-variance efficient *à la* Markowitz. It is shown that, when the market coefficients are deterministic functions of time, a mean-variance efficient portfolio realizes the (discounted) targeted return on or before the terminal date with a probability greater than 0.8072. This number is universal irrespective of the market parameters, the targeted return and the length of the investment horizon.

**1. Introduction.** In his seminal work, Markowitz [9] proposed the mean-variance portfolio selection model for a single investment period, where an agent seeks to minimize the risk of his investment, measured by the variance of his return, subject to a given mean return. (In Markowitz's original setup, the model is formulated as a *multi*-objective optimization problem, namely, to maximize the mean return and minimize the variance of the return. There are multiple solutions to this problem, leading to the so-called *efficient frontier*. Mathematically, each solution can be recovered by solving a *single*-objective optimization problem where the variance is to be minimized while the return is constrained at a given level.) The dynamic extension of the Markowitz model, especially in continuous time, has been studied extensively in recent years; see, for example, [2, 6, 7, 8, 11, 12]. (In particular,

---

Received October 2006; revised March 2006.

[1]Supported in part by NUS Grant R-146-000-059-112.

[2]Supported in part by RGC Earmarked Grants CUHK418605, 4175/03E, and the Croucher Senior Research Fellowship.

refer to [10] and [2] for elaborative discussions on the history of the mean-variance model.) In many of these works explicit, analytic forms of efficient portfolios have been obtained.

In spite of being awarded a Nobel prize in 1990, the mean-variance model has received criticisms since its inception, the main criticism being using variance as a risk measure. This, in turn, has subsequently led to many alternative models, such as those of semi-variance or shortfall, downside risk and lower partial moment. Another criticism is that the mean-variance model uses mathematical expectation, in contrast with models such as VaR incorporating probability. The argument is that a model with expectation, by its very definition, appears to work well only on average over a large number of different sample paths, which has little relevance with a real-world investor who would experience only *one* sample path realization over a fixed investment horizon.

This represents a typical dilemma in Mathematical Finance (or in any applied mathematics area for that matter): on one hand, one needs to establish models that are analytically or numerically tractable and, on the other hand, the models must be relevant to the real world. Fortunately, more often than not, model tractability and relevance co-exist nicely. For continuous-time portfolio selection, the mean-variance model is so mathematically simple and elegant that it generally admits closed-form solutions. On the other hand, we will show *analytically* in this paper that, although the model is being optimized in the average (expectation) sense, if one follows an efficient portfolio generated by the mean-variance model, then there is more than an 80% chance that he will reach his goal (the targeted return) on or before the prescribed terminal date. Moreover, this "goal-achieving" probability, 0.8072 to be more precise, is independent of the market parameters, the target or the length of the investment horizon. This astonishing 80% rule would provide reference and guidance in investment practice: one could simply follow a mean-variance efficient strategy, stop (i.e., withdraw from the stock market) as soon as his wealth hits the discounted value of the target (the chance of this happening is more than 80%); otherwise just follow the original strategy till the terminal time. This implied policy would meet the original target with a probability of more than 80%.

The 80% rule will be derived in this paper based on the known, explicit form of an efficient strategy, the probability distribution of the hitting time of an Itô process on a certain level, as well as some delicate optimization techniques.

It should be noted that the main technical assumption of the 80% rule is that the market coefficients (a.k.a. the investment opportunity set) are *deterministic* functions of time. It remains an interesting open question whether the result carries over to the case of stochastic coefficients and, if the answer is no, what the corresponding probability is.

The remainder of the paper is organized as follows. In Section 2 the continuous-time mean-variance model is formulated and its solution presented. Section 3 is devoted to proving the main result of the paper—the 80% rule. Some discussions and suggestions of possible open problems are given in Section 4.

**2. Mean-variance model and solution.**   Throughout this paper $(\Omega, \mathcal{F}, P, \{\mathcal{F}_t\}_{t \geq 0})$ is a fixed filtered complete probability space on which is defined a standard $\mathcal{F}_t$-adapted $m$-dimensional Brownian motion $\{W(t), t \geq 0\}$ with $W(t) \equiv (W^1(t), \ldots, W^m(t))'$ and $W(0) = 0$, and $T > 0$ is given and fixed representing the terminal time of an investment. In addition, we use $M'$ to denote the transpose of any vector or matrix $M$, and $L^2_{\mathcal{F}}(0, T; \Re^d)$ to denote the set of all $\Re^d$-valued, $\mathcal{F}_t$-progressively measurable stochastic processes $f(t)$ with $E \int_0^T |f(t)|^2 \, dt < +\infty$.

There is a capital market in which $m + 1$ *basic securities* (or *assets*) are traded continuously. One of the securities is a risk-free bank account whose value process $S_0(t)$ is subject to the following ordinary differential equation:

$$
\begin{aligned}
(2.1) \qquad & dS_0(t) = r(t)S_0(t) \, dt, \qquad t \geq 0, \\
& S_0(0) = s_0 > 0,
\end{aligned}
$$

where $r(t) > 0$ is the interest rate. The other $m$ assets are risky stocks whose price processes $S_1(t), \ldots, S_m(t)$ satisfy the following stochastic differential equation (SDE):

$$
\begin{aligned}
(2.2) \qquad & dS_i(t) = S_i(t)\left[ \mu_i(t) \, dt + \sum_{j=1}^m \sigma_{ij}(t) \, dW^j(t) \right], \qquad t \geq 0, \\
& S_i(0) = s_i > 0, \qquad i = 1, 2, \ldots, m,
\end{aligned}
$$

where $\mu_i(t)$ is the appreciation rate, and $\sigma_{ij}(t)$ is the volatility or dispersion rate of the stocks. We assume that all the given market parameters $r(t), \mu_i(t)$ and $\sigma_{ij}(t)$ are deterministic functions in $t \geq 0$.

Consider an agent, with an initial endowment $x_0 > 0$ and an investment horizon $[0, T]$, whose total wealth at time $t \in [0, T]$ is denoted by $x(t)$. Assume that the trading of shares is self-financed and takes place continuously, and that transaction cost and consumptions are not considered. Then $x(\cdot)$ satisfies (see, e.g., [4])

$$
dx(t) = \left\{ r(t)x(t) + \sum_{i=1}^m [\mu_i(t) - r(t)]\pi_i(t) \right\} dt
$$

(2.3) $$+ \sum_{j=1}^{m} \sum_{i=1}^{m} \sigma_{ij}(t)\pi_i(t)\,dW^j(t), \qquad 0 \le t \le T,$$

$$x(0) = x_0,$$

where $\pi_i(t), i = 1, 2, \ldots, m,$ denotes the total market value of the agent's wealth in the $i$th stock. We call the process $\pi(t) := (\pi_1(t), \ldots, \pi_m(t))', \ 0 \le t \le T,$ a *portfolio* of the agent.

DEFINITION 2.1.   A portfolio $\pi(\cdot)$ is said to be *admissible* if $\pi(\cdot) \in L^2_{\mathcal{F}}(0, T;$ $\Re^m)$ and the SDE (2.3) has a unique solution $x(\cdot)$ corresponding to $\pi(\cdot)$.

The agent's objective is to find an admissible portfolio $\pi(\cdot)$, among all admissible portfolios such that their expected terminal wealth $Ex(T) = z$, where $z \ge x_0 e^{\int_0^T r(t)\,dt}$ is given a priori, so that the risk measured by the variance of the terminal wealth

(2.4) $$\operatorname{Var} x(T) := E[x(T) - Ex(T)]^2 \equiv E[x(T) - z]^2$$

is minimized. The problem of finding such a portfolio $\pi(\cdot)$ is referred to as the *mean-variance portfolio selection problem*. Mathematically, we have the following formulation.

DEFINITION 2.2.   The mean-variance portfolio selection problem, with respect to the initial wealth $x_0$, is formulated as a constrained stochastic optimization problem parameterized by $z \ge x_0 e^{\int_0^T r(t)\,dt}$:

(2.5)
$$\text{minimize} \quad J_{\mathrm{MV}}(x_0; \pi(\cdot)) := E[x(T) - z]^2,$$
$$\text{subject to} \quad \begin{cases} x(0) = x_0, & Ex(T) = z, \\ (x(\cdot), \pi(\cdot)) & \text{admissible.} \end{cases}$$

The problem is called *feasible* (with respect to $z$) if there is at least one admissible portfolio satisfying $Ex(T) = z$. An optimal portfolio, if it ever exists, is called an *efficient portfolio* with respect to $z$.

REMARK 2.1.   In the formulation above the parameter $z$ is restricted to be no less than $x_0 e^{\int_0^T r(t)\,dt}$, which is the terminal payoff if all the initial wealth is put into the bank account. Hence, as standard with the single-period case, we are interested only in the nonsatiation portion of the minimum-variance set.

Define the covariance matrix $\sigma(t) := (\sigma_{ij}(t))_{m \times m}$. We impose the first basic assumption of this paper, which is essentially a uniform elliptic condition on the covariance matrix.

ASSUMPTION A1.   $\sigma(t)\sigma(t)' \geq \delta I \ \forall t \in [0,T]$ for some $\delta > 0$.

Next, we introduce the following notation:

$$(2.6) \qquad B(t) := (\mu_1(t) - r(t), \ldots, \mu_m(t) - r(t))$$

and

$$(2.7) \qquad \theta(t) \equiv (\theta_1(t), \ldots, \theta_m(t)) := B(t)(\sigma(t)')^{-1}.$$

ASSUMPTION A2.   $0 < \int_0^T |\theta(t)|^2 \, dt < +\infty$.

REMARK 2.2.   Assumption A1 and that $\int_0^T |\theta(t)|^2 \, dt < +\infty$ are to ensure that the market is arbitrage-free and complete, whereas that $\int_0^T |\theta(t)|^2 \, dt > 0$ is to guarantee that the mean-variance problem is feasible for *any* $z > 0$ (see [2], Theorem 3.1).

Assumptions A1 and A2 will be in force from now on.

With the above notation, equation (2.3) can be rewritten as

$$(2.8) \qquad \begin{aligned} dx(t) &= [r(t)x(t) + B(t)\pi(t)] \, dt + \pi(t)'\sigma(t) \, dW(t), \qquad 0 \leq t \leq T, \\ x(0) &= x_0. \end{aligned}$$

The following result, first derived in [12], gives a complete solution to the mean-variance portfolio selection problem.

THEOREM 2.1.   *The efficient portfolio corresponding to each given* $z \geq x_0 e^{\int_0^T r(t) \, dt}$ *can be uniquely represented as a feedback strategy*

$$(2.9) \qquad \begin{aligned} \pi^z(t) &\equiv (\pi_1^z(t), \ldots, \pi_m^z(t))' \\ &= -[\sigma(t)\sigma(t)']^{-1} B(t)'[x^z(t) - \gamma e^{-\int_t^T r(s) \, ds}], \end{aligned}$$

*where* $x^z(\cdot)$ *is the corresponding wealth process and*

$$(2.10) \qquad \gamma := \frac{z - x_0 e^{\int_0^T [r(t) - |\theta(t)|^2] \, dt}}{1 - e^{-\int_0^T |\theta(t)|^2 \, dt}} \geq z > 0.$$

*Moreover, the corresponding minimum variance can be expressed as*

$$(2.11) \qquad \operatorname{Var} x^z(T) = \frac{1}{e^{\int_0^T |\theta(t)|^2 \, dt} - 1}[z - x_0 e^{\int_0^T r(t) \, dt}]^2, \qquad z \geq x_0 e^{\int_0^T r(t) \, dt}.$$

**3. Probability of goal-achieving: The 80% rule.** In this section we answer the following question: given a target $z > x_0 e^{\int_0^T r(t)\,dt}$ (the one corresponding to $z = x_0 e^{\int_0^T r(t)\,dt}$ is the risk-free portfolio, hence, not interesting), if one follows the efficient strategy as stipulated in Theorem 2.1, what is the probability that the corresponding wealth reaches the discounted value of $z$ on or before $T$?

Let $x^z(\cdot)$ be the wealth process under the efficient portfolio corresponding to $z > x_0 e^{\int_0^T r(t)\,dt}$ as specified by Theorem 2.1. Define the first hitting time of the wealth on the discounted value of $z$:

$$(3.1) \qquad \tau^z := \inf\{0 \le t \le T : x^z(t) = z e^{-\int_t^T r(s)\,ds}\},$$

where (and throughout the paper) $\inf \varnothing := +\infty$.

THEOREM 3.1. *For any $z > x_0 e^{\int_0^T r(t)\,dt}$,*

$$(3.2) \quad \tau^z = \inf\left\{0 \le t \le T : \frac{3}{2}\int_0^t |\theta(s)|^2\,ds + \int_0^t \theta(s)\,dW(s) = \int_0^T |\theta(s)|^2 ds\right\}.$$

PROOF. Set $y(t) := x^z(t) - \gamma e^{-\int_t^T r(s)\,ds}$. Using the wealth equation (2.8) that $x^z(\cdot)$ satisfies and the fact that $\pi^z(t) = -[\sigma(t)\sigma(t)']^{-1}B(t)'y(t)$, we deduce

$$dy(t) = [r(t) - |\theta(t)|^2]y(t)\,dt - \theta(t)y(t)\,dW(t), \qquad 0 \le t \le T,$$

$$y(0) = \frac{x_0 - z e^{-\int_0^T r(t)\,dt}}{1 - e^{-\int_0^T |\theta(t)|^2\,dt}}.$$

The above equation has a unique solution

$$y(t) = y(0)\exp\left(\int_0^t [r(s) - \tfrac{3}{2}|\theta(s)|^2]\,ds - \int_0^t \theta(s)\,dW(s)\right), \qquad 0 \le t \le T.$$

Hence,

$$x^z(t) - z e^{-\int_t^T r(s)\,ds} = y(t) + (\gamma - z)e^{-\int_t^T r(s)\,ds}$$

$$= \frac{e^{-\int_t^T r(s)\,ds}(z - x_0 e^{\int_0^T r(s)\,ds})}{e^{\int_0^T |\theta(t)|^2\,dt} - 1}$$

$$\times [1 - e^{\int_0^T |\theta(t)|^2\,dt} e^{-(3/2)\int_0^t |\theta(s)|^2\,ds - \int_0^t \theta(s)\,dW(s)}].$$

Since $\dfrac{e^{-\int_t^T r(s)\,ds}(z - x_0 e^{\int_0^T r(s)\,ds})}{e^{\int_0^T |\theta(t)|^2\,dt} - 1} > 0$, we conclude that $x^z(t) - z e^{-\int_t^T r(s)\,ds} = 0$ if and only if the term in the above bracket vanishes, or

$$\tfrac{3}{2}\int_0^t |\theta(s)|^2\,ds + \int_0^t \theta(s)\,dW(s) = \int_0^T |\theta(s)|^2\,ds.$$

This proves (3.2).  □

It is interesting to note that the hitting time depends entirely on the behavior of the market as represented by $\theta(\cdot)$ (the market price of risk), and does *not* depend on the target $z$. The following result gives an analytical formula for calculating the probability that the hitting occurs on or before the terminal time.

Recall the error function of the standard normal distribution

$$(3.3) \qquad \mathrm{Erfc}(x) := \frac{2}{\sqrt{\pi}}\int_x^\infty e^{-v^2}\,dv, \qquad x \in \Re.$$

THEOREM 3.2.  *The probability that an efficient wealth process* $x^z(\cdot)$, *corresponding to* $z > x_0 e^{\int_0^T r(t)\,dt}$, *reaches the discounted value of* $z$ *on or before the terminal time* $T$ *is given by*

$$
\begin{aligned}
(3.4) \qquad P(\tau^z \le T) = {}& \frac{1}{2}\mathrm{Erfc}\left(-\frac{\sqrt{\int_0^T |\theta(s)|^2\,ds}}{2\sqrt{2}}\right) \\
& + \frac{1}{2}e^{3\int_0^T |\theta(s)|^2\,ds}\mathrm{Erfc}\left(\frac{5\sqrt{\int_0^T |\theta(s)|^2\,ds}}{2\sqrt{2}}\right).
\end{aligned}
$$

PROOF.  By Theorem 3.1,

$$\tau^z \equiv \tau = \inf\left\{0 \le t \le T : \varphi(t) = \int_0^T |\theta(s)|^2\,ds\right\},$$

where

$$\varphi(t) := \tfrac{3}{2}\int_0^t |\theta(s)|^2\,ds + \int_0^t \theta(s)\,dW(s), \qquad 0 \le t \le T.$$

By virtue of a time-change technique (see, e.g., [5]), there exists a one-dimensional standard Brownian motion $\widehat{W}(t)$, $t \ge 0$, on $(\Omega, \mathcal{F}, P)$ such that

$$\int_0^t \theta(s)\,dW(s) = \widehat{W}(\beta(t)), \qquad 0 \le t \le T,$$

where $\beta(t) := \int_0^t |\theta(s)|^2\,ds$. Hence,

$$\varphi(t) = \tfrac{3}{2}\beta(t) + \widehat{W}(\beta(t)), \qquad 0 \le t \le T.$$

Now,

$$P(\tau \leq T) = P\left(\sup_{0 \leq t \leq T} \varphi(t) \geq \beta(T)\right)$$

$$= P\left(\sup_{0 \leq t \leq \beta(T)} (\tfrac{3}{2}t + \widehat{W}(t)) \geq \beta(T)\right).$$

According to Borodin and Salminen [3], page 250, 1.1.4, the above probability equals

$$P(\tau \leq T) = \frac{1}{2}\operatorname{Erfc}\left(\frac{\beta(T)}{\sqrt{2\beta(T)}} - \frac{(3/2)\sqrt{\beta(T)}}{\sqrt{2}}\right)$$

$$+ \frac{1}{2}e^{3\beta(T)}\operatorname{Erfc}\left(\frac{\beta(T)}{\sqrt{2\beta(T)}} + \frac{(3/2)\sqrt{\beta(T)}}{\sqrt{2}}\right)$$

$$= \frac{1}{2}\operatorname{Erfc}\left(-\frac{\sqrt{\beta(T)}}{2\sqrt{2}}\right) + \frac{1}{2}e^{3\beta(T)}\operatorname{Erfc}\left(\frac{5\sqrt{\beta(T)}}{2\sqrt{2}}\right).$$

Since $\beta(T) = \int_0^T |\theta(s)|^2 \, ds$, the preceding expression is identical to (3.4). $\quad\square$

Define the following function:

$$(3.5) \qquad f(x) := \frac{1}{2}\operatorname{Erfc}\left(-\frac{x}{2\sqrt{2}}\right) + \frac{1}{2}e^{3x^2}\operatorname{Erfc}\left(\frac{5x}{2\sqrt{2}}\right), \qquad x \geq 0.$$

Theorem 3.2 states that

$$P(\tau^z \leq T) = f\left(\sqrt{\int_0^T |\theta(s)|^2 \, ds}\right).$$

We plot $f$ in Figure 1. By inspection, $f$ has a minimum value slightly above 0.80. We now *prove* this *analytically*.

Lemma 3.1. *The function $f$ defined by (3.5) satisfies*

$$(3.6) \qquad f(x) \geq N\left(\frac{1}{\sqrt{5}}\right) + \frac{1}{12}\sqrt{\frac{10}{\pi}}e^{-1/10} \approx 0.8072 \qquad \forall x \geq 0.$$

Proof. Making use of the relation $\operatorname{Erfc}(x) \equiv 2(1 - N(\sqrt{2}x))$, where $N(x) := \frac{1}{\sqrt{2\pi}}\int_{-\infty}^x e^{-v^2/2} \, dv$ is the p.d.f. of the standard normal distribution, we rewrite $f$ as

$$(3.7) \qquad f(x) = N\left(\frac{x}{2}\right) + e^{3x^2}\left(1 - N\left(\frac{5x}{2}\right)\right), \qquad x \geq 0.$$

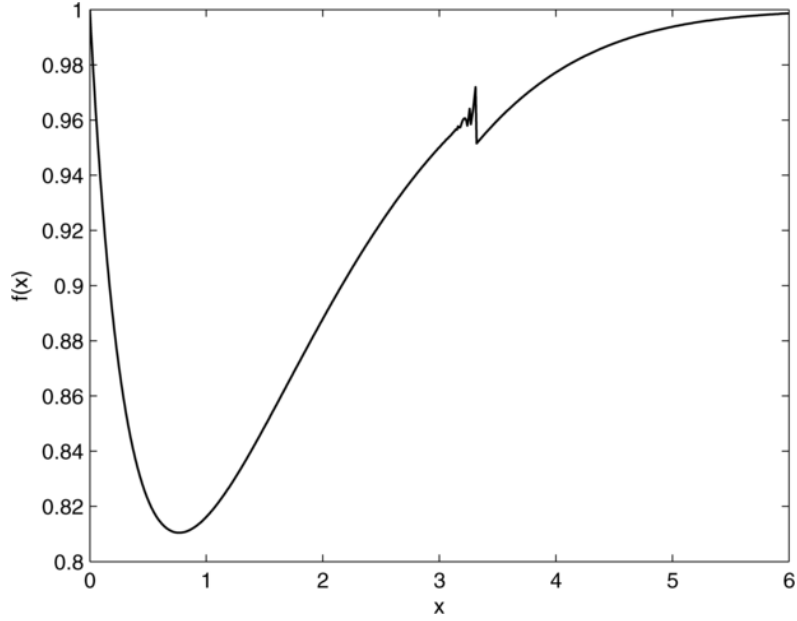

Fig. 1.

We first prove that

$$(3.8) \qquad f(x) \geq N\left(\frac{1}{2}\right) + \frac{1}{\sqrt{2\pi}} \frac{\sqrt{41}-5}{4} e^{-1/8} \approx 0.8150 \qquad \forall x \geq 1.$$

Indeed, employing the following estimate for the function $N(\cdot)$ (see page 933, 26.2.24, of [1]),

$$N(x) \leq 1 - \frac{\sqrt{4+x^2}-x}{2\sqrt{2\pi}} e^{-x^2/2} \qquad \forall x > 1.4,$$

we have

$$(3.9) \qquad \begin{aligned} f(x) &\geq N\left(\frac{x}{2}\right) + e^{3x^2} \frac{\sqrt{16+25x^2}-5x}{4\sqrt{2\pi}} e^{-25x^2/8} \\ &= N\left(\frac{x}{2}\right) + \frac{\sqrt{16+25x^2}-5x}{4\sqrt{2\pi}} e^{-x^2/8} \qquad \forall x \geq 1. \end{aligned}$$

Denote

$$g(x) := N\left(\frac{x}{2}\right) + \frac{\sqrt{16+25x^2}-5x}{4\sqrt{2\pi}} e^{-x^2/8}, \qquad x \geq 0.$$

Then its derivative is (after some manipulations)

$$\dot{g}(x) = \frac{e^{-x^2/8}}{\sqrt{2\pi}} \left[ \left( \frac{5x^2}{16} - \frac{25x^3}{16\sqrt{16+25x^2}} \right) + \left( \frac{21x}{4\sqrt{16+25x^2}} - \frac{3}{4} \right) \right].$$

Now,

$$\frac{5x^2}{16} - \frac{25x^3}{16\sqrt{16+25x^2}} \geq \frac{5x^2}{16} - \frac{25x^3}{16\sqrt{25x^2}} = 0 \qquad \forall x \geq 0.$$

On the other hand, since $\frac{x}{\sqrt{16+25x^2}}$ is strictly increasing in $x \geq 1$,

$$\frac{21x}{4\sqrt{16+25x^2}} - \frac{3}{4} \geq \frac{21}{4\sqrt{41}} - \frac{3}{4} > 0 \qquad \forall x \geq 1.$$

Therefore, $\dot{g}(x) > 0$ for $x \geq 1$, that is, $g$ is strictly increasing in $x \geq 1$. This implies (3.8) taking note of (3.9).

Next, consider the continuous function $f$, which must admit a minimum point on $[0,1]$. The candidates for such a minimum point are $0$, $1$ or $x^* \in (0,1)$ satisfying $\dot{f}(x^*) = 0$, or

$$(3.10) \qquad e^{3x^{*2}}\left(1 - N\left(\frac{5x^*}{2}\right)\right) = \frac{1}{3\sqrt{2\pi}x^*}e^{-x^{*2}/8}.$$

Note that such an $x^*$ does not need to exist on $(0,1)$; but if indeed it exists, then necessarily

$$f(x^*) \equiv N\left(\frac{x^*}{2}\right) + e^{3x^{*2}}\left(1 - N\left(\frac{5x^*}{2}\right)\right)$$

$$= N\left(\frac{x^*}{2}\right) + \frac{1}{3\sqrt{2\pi}x^*}e^{-x^{*2}/8}.$$

To estimate the above value, define

$$h(x) := N\left(\frac{x}{2}\right) + \frac{1}{3\sqrt{2\pi}x}e^{-x^2/8}, \qquad x > 0.$$

Then $\dot{h}(x)$ has a unique root $\hat{x} = \frac{2}{\sqrt{5}}$. Moreover, $\dot{h}(x) < 0$ for $0 < x < \hat{x}$ and $\dot{h}(x) > 0$ for $x > \hat{x}$. Hence, $\hat{x}$ must be the global minimum of $h$, which implies

$$(3.11) \quad h(x) \geq h\left(\frac{2}{\sqrt{5}}\right) = N\left(\frac{1}{\sqrt{5}}\right) + \frac{1}{12}\sqrt{\frac{10}{\pi}}e^{-1/10} \approx 0.8072 \qquad \forall x > 0.$$

As a result,

$$f(x^*) \equiv h(x^*) \geq 0.8072.$$

However, $f(0) = 1$ and $f(1) \approx 0.8162$; so we conclude that the minimum value of $f$ on $[0,1]$ is at least $0.8072$. By virtue of (3.8), we arrive at

$$f(x) \geq 0.8072 \qquad \forall x \geq 0. \qquad\qquad \square$$

THEOREM 3.3. *We have the following lower bound:*

$$(3.12) \qquad P(\tau^z \leq T) \geq N\left(\frac{1}{\sqrt{5}}\right) + \frac{1}{12}\sqrt{\frac{10}{\pi}}e^{-1/10} \approx 0.8072.$$

PROOF.   This is immediate from Theorem 3.2 and Lemma 3.1.   □

The above lower bound, $N(\frac{1}{\sqrt{5}}) + \frac{1}{12}\sqrt{\frac{10}{\pi}}e^{-1/10} \approx 0.8072$, has been obtained *analytically*. Notice that it is not necessarily the tightest bound, as it was derived based on the global minimum of $h$ (which was obtained analytically), rather than that of $f$ (which seems to be impossible to get analytically). However, 0.8072 is already a very good lower bound because it is only slightly smaller than the minimum value of $f$ as suggested by Figure 1.

**4. Discussions.**   The main results derived in Section 3 are quite surprising and counter-intuitive. Mean-variance portfolio selection model, as many other stochastic optimization models, in its nature is one based on averaging over all the possible random scenarios; so an optimal solution is optimal only in the sense of an average. There have been debates on the sensibility of the model in terms of how much its solutions could guide real investment in practice. Moreover, due to the presence of the variance in its objective, the mean-variance model is *not* compatible with the dynamic programming principle. (The dynamic programming principle holds if we remove the constraint on the mean by introducing a suitable $\lambda$; nevertheless, the new problem is an auxiliary problem which is *not* completely equivalent to the original problem (although the latter can be solved based on the solution to the former; see [6] and [12]).) Specifically, an optimal portfolio generated initially may no longer be optimal half way through. Thus, one may tend to consider the mean-variance model to be "unfavorable," particularly in the dynamic setting, in not being able to generate sound investment policies. Now, the 80% rule demonstrates that there is a very high chance that a mean-variance strategy would lead to the full realization of the prescribed financial goal. Hence, while mean-variance arguably may not be the best model for portfolio management, it could indeed generate sound solutions.

The analytically derived *lower bound* of the probability, approximately 0.8072, is universal independent of the target $z$, the market represented by $\theta(\cdot)$ and the time horizon $T$. Theorem 3.2 asserts that the *exact* probability does depend on $\int_0^T |\theta(s)|^2 \, ds$ as an aggregation of the market and the horizon. Therefore, given a market, one could carefully choose the investment horizon $T$ so as to further increase the "goal-achieving probability" (after all, how long an investor is going to plan his investment is a *part* of his overall decision). Some insight in this aspect could be obtained from the analytical formula (3.4) together with its graph in Figure 1.

Having said all these, one should bear in mind that the results obtained in this paper are *theoretical* ones based on a number of assumptions, including that the stock prices are driven by Brownian motions, the market is complete, the market coefficients are deterministic and the transaction costs are ignored. Some of these assumptions may be purely technical and some may be essential. This, however, should not be a concern when the results are interpreted with care and discretion, and the model is used properly as a reference or study tool (as opposed to a trading tool) in the same spirit as the Black–Scholes model for option pricing. In fact, from another angle, this 80% rule could serve as a test on the validity of the aforementioned hypotheses: if the rule fails in, say, an extensive empirical study, then it might be an indication that one or more of the hypotheses are rejected by the data. In summary, we feel that the results, including the 80% rule, can shed some light on the portfolio theory and offer guidance and reference for investment practice.

We conclude this paper by pointing out some interesting open questions. First of all, the model in our paper allows for bankruptcy, that is, the wealth is allowed to go negative (e.g., borrow from bank and continue trading by buying stock on margin). Although the 80% rule dictates that in the end there is a high chance to reach one's goal, there could be an equally high chance that one has to experience bankruptcy *before* the goal is ever reached. It is easy to calculate the probability of going bankrupt before the terminal time, as well as the probability that this will happen before the goal is reached. The latter would depend on the length of the investment horizon, $T$; hence, one may determine $T$ so as to minimize the probability that bankruptcy occurs earlier than the goal-achieving. This would be an interesting problem. Nonetheless, an even more interesting problem is to consider the class of admissible portfolios which exclude bankruptcy in the first place. The mean-variance portfolio selection with bankruptcy prohibition has been solved in [2] in a very general setting, and explicit solutions have been obtained for the case of deterministic market coefficients. An open problem is therefore what the goal-achieving probability is for such a no-bankruptcy efficient portfolio.

Yet another, perhaps more challenging, open problem is to consider the case where all the market coefficients are stochastic. While the corresponding efficient portfolios have been obtained in [8] and [7], the estimation of hitting probability would require more delicate stochastic analysis.

**Acknowledgments.** We are grateful to Duan Li for suggesting incorporating a stopping rule into an efficient strategy, which motivated this research, and to Hanqing Jin for useful discussions. We also thank the participants in the various conferences at which previous versions of this paper were presented. Last but not least, we thank an Associate Editor and the two anonymous referees for constructive comments that have led to a much improved version.

DEPARTMENT OF MATHEMATICS
NATIONAL UNIVERSITY OF SINGAPORE
2 SCIENCE DRIVE 2
SINGAPORE 117543
REPUBLIC OF SINGAPORE
E-MAIL: matlx@nus.edu.sg

DEPARTMENT OF SYSTEMS ENGINEERING
  AND ENGINEERING MANAGEMENT
CHINESE UNIVERSITY OF HONG KONG
SHATIN
HONG KONG
E-MAIL: xyzhou@se.cuhk.edu.hk
URL: http://www.se.cuhk.edu.hk/~xyzhou